\documentclass[letterpaper,12pt]{article}
\usepackage{amssymb,latexsym}
\usepackage{enumerate}
\usepackage{amsthm}
\usepackage{amsmath,amsfonts}
\usepackage{graphicx}
\usepackage{xcolor,graphicx}
\usepackage{epsfig,bm,epsf,float}
\usepackage{footnote}
\theoremstyle{plain}

\newtheorem*{xthm}{Theorem}

\newtheorem*{xconj}{Conjecture}

\theoremstyle{definition}

\theoremstyle{remark}

\numberwithin{equation}{section}
\begin{document}
\title{Some History about Twin Prime Conjecture}
\maketitle
\section{Introduction}
This article is a collected information from some books and papers, and in most cases the original sentences is reserved.
\subsection{Historical Problems}
There exist infinitely many primes that do not belong to any pair of twin primes.
[Hint: Consider the arithmetic progression 2lk + 5 for k = 1,2,....]\cite[p. 59]{Burton}

The twin primes 5 and 7 are such that one half their sum is a perfect number. Are there
any other twin primes with this property?
[Hint: Given the twin primes p and p + 2, with $p > 5$, 1/2(p + p + 2) = 6k for some
$k > 1$.]\cite[p. 225]{Burton}

There are many questions concerning primes which have resisted all assaults for two centuries and more.\cite[p. 6]{Leveque}
лллллллллллллллллллллл
The set {(3, 5), (5, 7), (11, 13), ... } of twin prime pairs (q, q+ 2) has been studied
by Brun (1919), Hardy and Littlewood (1922), Selmer (1942), Fr\"{o}berg (1961),
Weintraub (1973), Bohman (1973), Shanks and Wrench (1974), and Brent (1975,
1976). Currently M. Kutrib and J. Richstein (1995) are completing a study similar
to the present one.\cite{Nicely}
лллллллллллллллл

Two of the oldest problems in the theory of numbers, and indeed in the whole of mathematics, are the so-called \emph{twin prime conjecture} and \emph{binary Goldbach problem}. The first of these asserts that

\emph{``There exist infinitely many primes $p$ such that $p + 2$ is a prim''}

and the second that

\emph{``every even natural number $N$ from some point onward can be expressed as the sum of two primes.''} \cite[p. 1]{Halberstam}\\
лллллллллллллллл

One charm of the integers is that easily stated problems, which often sound simple, are often very difficult and sometimes even hopeless given the state of our current knowledge. For instance, in a 1912 lecture at an international mathematical congress, Edmund Landau mentioned four old conjectures that appeared hopeless at that time:\\
\begin{itemize}
  \item Every even number greater than two is the sum of two primes.
  \item Between any consecutive squares there is a prime number.
  \item There are infinitely many twin primes.
  \item There are infinitely many primes of the form $n^2 + 1$.
\end{itemize}

Today, his statement could be reiterated. During the more than 80 years that have passed, much intensive research has been conducted on all of these conjectures, and we now know that every large enough odd number is the sum of three primes, and every even number can be written as the sum of a prime and a number that is divisible by at most two different primes; there are infinitely many consecutive odd integers for which one is a prime and the other has at most two prime divisors; for infinitely many n, there is a prime between $n^2$ and $n^2+n^{1.1}$; for infinitely many n, $n^2+ 1$ has at most two distinct prime divisors. Unfortunately, the methods used to achieve these rather deep
results cannot be generalized to prove the more general conjectures, and given the state of mathematics today, the resolution of these conjectures can still be called hopeless.\cite[pp. v, vi]{Erdos1}
лллллллллллллллл

The conjecture that $n^2 + 1$ is prime for infinitely many integers $n$ is sometimes called Landau's conjecture, since he discussed it in his address to the 1912 International Congress as one of the particularly challenging unsolved problems about prime numbers (the other three being the Goldbach conjecture, the infinitude of twin primes, and the existence of a prime between any two consecutive squares), cf. Proc., Fifth Int. Cong, of Math., vol. 1, Cambridge, 1913, pp. 93-108, particularly p. 106; also in \cite{Landau c}, vol. 5, pp. 240-255, particularly p. 253, and in Jahresber. Deutsch. Math. Verein., vol. 21 (1912), pp. 208-228, particularly p. 224. While these four problems remain unsolved as stated, interesting partial results have been obtained for each of them. The nearest approach to Landau's conjecture to date is H. Iwaniec's result (\cite{Iwaniec}): that there are infinitely many integers $n$ such that $n^2+1$  is the product of at most two primes.\cite[p. 337]{Bateman}
лллллллллллллллл

Although there is an infinitude of primes, their distribution within the positive integers is most mystifying. Repeatedly in their distribution we find hints or, as it were, shadows of a pattern; yet an actual pattern amenable to precise description remains elusive. The difference between consecutive primes can be small, as with the pairs 11 and 13, 17 and 19, or for that matter 1,000,000,000,061 and 1,000,000,000,063. At the same time there exist arbitrarily long intervals in the sequence of integers that are totally devoid of any primes.

Consecutive primes cannot only be close together, but also can be far apart; that is, arbitrarily large gaps can occur between consecutive primes. Stated precisely: Given any positive integer n, there exist n consecutive integers, all of which are composite. To prove this, we simply need to consider the integers
$$
(n+1)!+2, (n+1)!+3,\ldots,(n+1)!+(n+1)
$$
For instance, if a sequence of four consecutive composite integers is desired, then the previous argument produces 122, 123, 124, and 125. Of course, we can find other sets of four consecutive composites, such as 24, 25, 26,
27 or 32, 33, 34, 35.

As this example suggests, our procedure for constructing gaps between two consecutive primes gives a gross overestimate of where they occur among the integers. The first occurrences of prime gaps of specific lengths, where all the intervening integers are composite, have been the subject of computer searches. For instance, there is a gap of length 778 (that is, $p_{n+1} - p_n = 778$) following the prime 42,842,283,925,351. No gap of this size exists between two smaller primes. The largest effectively calculated gap between consecutive prime numbers has length 1132, with a string of 1131 composites immediately after the prime 1,693,182,318,746,371.\\
Interestingly, computer researchers have not identified gaps of every possible width up to 1132. The smallest missing gap size is 796. The conjecture is that there is a prime gap (a string of 2k-1 consecutive composites between two primes) for every even integer 2k.

This brings us to another unsolved problem concerning the primes, the Goldbach conjecture. In a letter to Leonhard Euler in the year 1742, Christian Goldbach hazarded the guess that every even integer is the sum of two numbers that are either primes or 1. A somewhat more general formulation is that every even integer greater than 4 can be written as a sum of two odd prime numbers. This is easy to confirm for the first few even integers:
2=1+1\\
4=2+2=1+3\\
6=3+3=1+5\\
8=3+5=1+7\\
10 = 3 + 7 = 5 + 5\\
12 = 5 + 7= 1 + 11\\
14 = 3 + 11 =7 + 7= 1 + 13\\
16 = 3 + 13 = 5 + 11\\
18 = 5 + 13 = 7 + 11 = 1 + 17\\
20 = 3 + 17 = 7 + 13 = 1 + 19\\
22 = 3 + 19 = 5 + 17= 11 + 11\\
24 = 5 + 19 = 7+17= 11 + 13 = 1+23\\
26 = 3 + 23 = 7+19= 13 + 13\\
28 = 5 + 23 = 11 + 17\\
30 = 7 + 23 = 11 + 19= 13 + 17= 1+29\\
Although it seems that Euler never tried to prove the result, upon writing to Goldbach at a later date, Euler countered with a conjecture of his own: Any even integer ($\geq6$) of the form 4n + 2 is a sum of two numbers each being either a prime of the form 4n + 1 or 1.

The numerical data suggesting the truth of Goldbach's conjecture are overwhelming. It has been verified by computers for all even integers less than $4\cdot10^{14}$. \cite[pp. 50, 51]{Burton}

As the integers become larger, the number of different ways in which 2n can be expressed as the sum of two primes increases. For example, there are 219,400 such representations for the even integer 100,000,000. Although this supports the feeling that Goldbach was correct in his conjecture, it is far from a mathematical proof, and all attempts to obtain a proof have been completely unsuccessful. One of the most famous number theorists of the last century, G. H. Hardy, in his address to the Mathematical Society of Copenhagen in 1921, stated that the Goldbach conjecture appeared ``... probably as difficult as any of the unsolved problems in mathematics.''

We remark that if the conjecture of Goldbach is true, then each odd number larger than 7 must be the sum of three odd primes. To see this, take n to be an odd integer greater than 7, so that n-3 is even and greater than 4; if n-3 could be expressed as the sum of two odd primes, then n would be the sum of three.

The first real progress on the conjecture in nearly 200 years was made by Hardy and Littlewood in 1922. On the basis of a certain unproved hypothesis, the so-called generalized Riemann hypothesis, they showed that every sufficiently large odd number is the sum of three odd primes. In 1937, the Russian mathematician I. M. Vinogradov was able to remove the dependence on the generalized Riemann hypothesis, thereby giving an unconditional proof of this result; that is to say, he established that all odd integers greater than some effectively computable $n_0$ can be written as the sum of three odd primes.

n = p1 + p2 + P3 (n odd, n sufficiently large)

Vinogradov was unable to decide how large $n_0$ should be, but Borozdkin (1956) proved that $n_0 < 3^{3^{15}}$. In 2002, the bound on $n_0$ was reduced to $10^{1346}$. It follows immediately that every even integer from some point on is the sum of either two or four primes. Thus, it is enough to answer the question for every odd integer n in the range $9\leq n\leq n_0$, which, for a given integer, becomes a matter of tedious computation (unfortunately, $n_0$ is so large that this exceeds the capabilities of the most modern electronic computers).

Because of the strong evidence in favor of the famous Goldbach conjecture, we readily become convinced that it is true. Nevertheless, it might be false. Vinogradov showed that if A(x) is the number of even integers $n\leq x$ that are not the sum of two primes, then
$$
\lim_{x\rightarrow\infty}\frac{A(x)}{x}=0
$$
This allows us to say that "almost all" even integers satisfy the conjecture. As Edmund Landau so aptly put it, "The Goldbach conjecture is false for at most $0\%$ of all even integers; this at most $0\%$ does not exclude, of course, the possibility that there are infinitely many exceptions."
\cite[p. 52]{Burton}
лллллллллллллллл

Goldbach problem

Every sufficiently large even number 2n is the sum of two primes. The asymptotic formula for the number of representations is
$$
r_2(2n)\sim C_2\frac{4n}{(\log2n)^2}\prod_{\substack{p>2\\p|n}}\frac{p-1}{p-2}
$$
with
$$
r_2(2n)=\#\{(p, q):\ p, q\ \text{primes},\ p+q=2n\}
$$
$$
C_2=\prod_{p>2}\left(1-\frac{1}{(p-1)^2}\right)=0.66016\ldots
$$
\cite[p. 404]{Ribenboim}
ллллллллллллллллллл

лллллллллллллллл
Recall that $(p, p+2)$ is a twin prime pair (or $p, p+2$ are twin primes) if $p$ and $p+2$ are both prime. Although the existence of such primes would have been easily understood by the ancient Greeks, there is no evidence that they were considered by them, or indeed by any mathematician, until the nineteenth century. The first mention of twin primes in the literature appears in de Polignac's paper of 1849, in which he speculates about the distribution of primes. \cite[p. 7]{Klyve}
лллллллллллллллллл

There is only one pair of prime numbers that differ by 1-namely, the primes 2 and 3- since all of the succeeding primes are odd. No matter how far out we go in a table of primes, we find pairs of primes that differ by 2, for example,
$$
3,5; 5,7; 11,13; 17,19; 29,31; \ldots ; 101,103; \ldots
$$
(Incidentally, each of these pairs, aside from the second, begins with a number greater than the second number of the preceding pair; this is because one of the numbers n, n+2, and n+4 is always divisible by 3, so that if these are all to be primes, we must have n=3.)

The question now arises whether there are infinitely many "prime pairs," that is, whether there are infinitely many numbers n for which n and n+2 are prime. The methods of number theory and of analysis have, to this day, not proven powerful enough to answer this question. (One would certainly place one's bet on a yes answer.)\cite[p. 94]{Landau}

\subsection{Sieve Methods}
Sieve theory has been used in attempts to prove that there exists infinitely many twin primes. Many authors have worked with this method.\cite[p. 264]{Ribenboim}\\[-7pt]

лллллллллллллллллл
The theory of sieves has become an important and sophisticated tool in the theory of numbers during the past 52 years.
Many important and deep results have been proved related to such unanswered questions as the Twin Primes Problem and Goldbach's Conjecture. For the most part sieves have been employed in problems concerning primes.\cite{Andrews}
лллллллллллллллллллл

Any arithmetical sieve (in our sense of the word) is based on the following idea: given a finite integer sequence $\mathfrak{A}$, a sequence $\mathfrak{B}$ of primes, and a number $z (\geq 2)$, if we remove (sift out) from $\mathfrak{A}$ all those elements that are divisible by primes of $\mathfrak{B}$ less than $z$, then the residual (unsifted) members of $\mathfrak{A}$ can have prime divisors from $\mathfrak{B}$ only if these are large (in the sense of being $\geq z$); a fortiori, each unsifted element of $\mathfrak{A}$ has only few such prime divisors
(provided that $z$ is not too small compared with $\max_{a\in\mathfrak{A}}|a|$). The object of any sieve theory is to estimate the number $S(\mathfrak{A}; \mathfrak{B}, z)$ of unsifted elements of $\mathfrak{A}$ or of some weighted form of this number.

There are many arithmetical investigations in which it suffices to have a good upper bound for this numberЧthe famous Brun-Titchmarsh inequality is a classic example. Other investigations require an asymptotic formula for $S(\mathfrak{A}; \mathfrak{B}, z)$, but in general this is an extremely difficult if not hopeless problem.

The first to devise an effective sieve method that goes substantially beyond the sieve of Eratosthenes was Viggo Brun. Brun's sieve is heavily combinatorial in character (most sieves share this property, at any rate to some degree); its complicated structure and Brun's own early accounts tended to discourage closer study and it is fair to say that only a small number of experts ever mastered its technique or realized its full potential. Hardy and Littlewood used Brun's
sieve only once, to derive an early version of the Brun-Titchmarsh inequality; and Hardy expressed the opinion that Brun's sieve did not seem ``..., sufficiently powerful or sufficiently profound to lead to a solution [of Goldbach's Problem]''.

Selberg published an account of his upper bound sieve in a short paper in 1947; and subsequently he indicated its place in the construction of lower bound sieves in several important expository articles, without actually ever publishing any details.\cite[p. 5, 6]{Halberstam}

However extensive a table of primes, prime twins or prime triplets\footnote{(n,n+2,n+6) are primes} may be, we cannot infer from such a table that there exist infinitely many terms in any one of these three sequences, however plausible each proposition might seem from the frequent occurrence of primes, twins or triplets in the table itself. The foregoing remarks indicate that, while the sieve of Eratosthenes is a useful method from a numerical point of view, it does not automatically yield also theoretical information.

It is clear how we might hope to gain such theoretical information from the sieve: if we sift the integer sequence $\mathfrak{A}$ (for example, think of $\mathfrak{A}$ as one of the sequences mentioned above) by the set of primes, $\mathfrak{B}$, we might be able to devise a means of estimating the number of elements of $\mathfrak{A}$ surviving the sifting procedure; thus an upper bound for this number would show that there cannot be too many primes (or prime twins, or prime triplets)Чinformation that is often useful in arithmetical investigationsЧor, better still, a positive lower bound (as $z\rightarrow\infty$) would show that there exist infinitely many primes (a fact easily proved in other ways), or infinitely many prime twins (one of the notorious unsolved problems of prime number theory).

What we require, therefore, is a theoretical analogue of the Eratosthenes sieve. Evidently certain important questions in prime number theory can be formulated in terms of sifting a certain type of sequence $\mathfrak{A}$ by a set $\mathfrak{B}$ of primes - the reader should have no difficulty now in expressing Goldbach's problem, or the problem of the infinitude of primes of the form $n^2 + 1$, as a sieve problem - and their answer depends on gaining significant information about the number of elements of $\mathfrak{A}$ surviving the sifting procedure. Indeed, there are many such questions; moreover, for many of these a sieve approach appears to offer the most promising, and sometimes the only, mode of attack. The construction of a theoretical sieve is therefore of great interest, and such interest is heightened by the fact that the construction is, beyond a certain stage of effectiveness, exceedingly difficult. The modern sieve is, as yet, an imperfect method, and one may doubt whether it alone will ever succeed in settling difficult questions such as the Goldbach and prime twins conjectures. It is noticed that even now deep results from prime number theory can be harnessed effectively to existing sieve methods; on the other band, the very generality of a sophisticated sieve theorem tends to militate against its success in a specific problem.\cite[p. 13]{Halberstam}

Brun attributes the term "prime twins" to St\"{a}ckel, who carried out some numerical calculations in connection with this and related problems. There are 1224 prime twins below 100,000 (Glaisher, Mess. Math. 8 (1878), 28-33) and 5328 below 600,000 (see Hardy-Littlewood, \emph{Note on Messrs. Shah and Wilson's paper entitled On an empirical formula connected with Goldbach's Theorem}, Proc. Cambridge Philos. Soc. 19 (1919), 245-254.). Schinzel and Sierpinski(\emph{Sur certaines hypoth\`{e}ses concernant les nombres premiers.} Acta Arith. 4 (1958), 185-208; Corrigendum: \emph{ibid}. 5 (1959), 259. MR 21,4936.) quote tables of prime``quadruplet'' and ``quintuplets'' from V. A. Golubev, \emph{Anzeigen Oesterr, Akad. Wiss.} (1956), 153-157; \emph{ibid} (1957), 82-87 and Sexton, \emph{ibid} (1955), 236-239. D. H. and E. Lehmer have carried out counts of various prime pairs, triplets and quadruplets up to $4\cdot10^7$, and their results are deposited in the \emph{Unpublished Math. Tables file of Math.} tables and other aids to computation. Shen, \emph{Nord. Tidskr. Inf. Beh.} 4 (1964), 243-245 reports that he has checked Goldbach's binary conjecture up to 33,000,000.

For earlier formulations of the general sieve problem see especially Ankeny and Onishi,  and the important recent memoir of Selberg. Of course, there are other accounts of sieve methods in which various degrees of generality and precision of formulation are achieved.\cite[p. 33]{Halberstam}

The prime twins conjecture is best known and simplest instance of the question where we ask whether a given irreducible polynomial with a prime argument assumes infinitely many prime values.\cite[p. 261]{Halberstam}

ллллллллллллллллллллл

\subsection{Brun's Constant}
The first decades of the twentieth century saw the establishment of major conjectures and theorems that today underlie our understanding of twin primes, along with a significant extension of the computation and enumeration of these numbers. Most important among these was the work of Brun, who provided the first non-trivial bound on twin primes, and the establishment of the conjecture of Hardy and Littlewood, given in their ground-breaking paper in 1923.

In 1919, Brun adapted and improved earlier work of J. Merlin(\emph{Sur quelque th\'{e}or\`{e}mes d'arithm\'{e}tique et un \'{e}nonc\'{e} qui les contient.} Comptes Rendus Acad. Sci. Paris, 153 (1911), pp. 516-518.) on the sieve of Eratosthenes to find the first non-trivial result concerning twin primes. Brun showed that
$$
\pi_2(x)=O\left(\frac{x(\log\log x)^2}{\log^2x}\right)
$$

In fact, Brun obtained an effective version of this bound. In particular, he showed that for some $x_0$ and all $x>x_0$,
$$
\pi_2(x)<7200\frac{x}{\log^2x}(\log\log x)^2+\frac{x}{\log^6x}+x^{3/4}
$$
Brun immediately followed this work with the announcement of a stronger bound(\emph{Le crible d'Eratosth\`{e}ne et le th\'{e}or\`{e}me de Goldbach.} Comptes Rendus Acad. Sci. Paris, 168 (1919), pp. 544-546.), in which he showed
$$
\pi_2(x)=O\left(\frac{x}{\log^2x}\right)
$$
and once again he found an effective version of his bound with
$$
\pi_2(x)<100\frac{x}{\log^2x},\qquad\text{for some}\  x > x_0,
$$
where $x_0$ is an effective computable constant(\emph{Le crible d'Eratosth\`{e}ne et le th\'{e}or\`{e}me de Goldbach.} Videnselsk. Skr. 1, No. 3 (1920)). BrunТs work was particularly important since even his first result above implied that the sum of the reciprocals of the twin primes converges.

Ever since BrunТs work, there has been a steady improvement in our knowledge of the function $\pi_2(x)$. Most of this work has taken the form of effective and ineffective arguments which seek to improve upon Brun. A word of explanation about this description is warranted:
\begin{itemize}
  \item A bound with numerical constants is \emph{explicit} if it holds at a certain given point; e.g.,
  $$
  \pi_2(x) < 16\alpha x/ \log^2x,\quad\text{ for all}\ x > 100.
  $$
  \item A bound is \emph{effective} if it holds starting at a given point, and if both the starting point and possibly the implied constants in the bound could be calculated in principle using the methods of the proof of the theorem.
  \item A bound is \emph{ineffective} if we don't know, even in principle, how to make it effective.
\end{itemize}
Because the conjectured asymptotic value of $\pi_2(x)$ is $2\alpha li_2(x)$, we usually express bounds as multiples of this value. That is, researchers seek the smallest $c$ such that $\pi_2(x)$ can be shown to be bounded by $c.2\alpha.li_2(x)$ for sufficiently large $x$. Brun himself obtained $c = 100/(2\alpha)\approx75.7$. Progress in obtaining ever-smaller values of $c$ can be seen in the following table, some of which was taken from Narkiewicz.\cite[pp. 13, 14]{Klyve}
\begin{table}[h!]
  \centering
  \caption{Estimating twin prime Constant}\label{constant}
\begin{tabular}{l l l}\\\\
Year &$c$& Name\\
\hline
1919&75.7\ldots&Brun\\
1947&8& Selberg\\
1964&6& Pan\\
1974&4& Halberstam, Richert\\
1978&3.9171& Chen\\
1983&34/9 = 3.777\ldots& Fouvry, Iwaniec\\
1984&64/17 = 3.764\ldots& Fouvry\\
1986 &3.5& Bombieri, Friedlander, Iwaniec\\
1986 &3.454\ldots &Fouvry, Grupp\\
1990 &3.418& Wu\\
2003 &3.406& Cai, Lu\\
2004 &3.3996& Wu
\end{tabular}
\end{table}

Each of these bounds holds for sufficiently large $x$, but for none of them has the starting point been made explicit. Explicit bounds are much harder to come by. In fact, we know of only one explicit bound on $\pi_2(x)$. This was given by Riesel and Vaughan as
$$
\pi_2(x)<\frac{16\alpha x}{(7.5+\log x)\log x},\qquad x>e^{42}
$$

\begin{table}[h!]
  \centering
  \caption{Estimating BrunТs Constant}\label{constant}
\begin{tabular}{l l l l}\\\\
Year &Computed to& Estimate of $B$& Author\\
\hline
1942 &$2\times10^5$    &$1.901 \pm0.0014$ &Selmer\\
1961 &$2^{20}$      &$1.90195 \pm 3 \times 10^{-5}$& Fr\"{o}berg\\
1973 &$2\times10^{9}$  &$1.90216 \pm 5 \times 10^{-6}$ &Bowman\\
1974 &$32,452,843$    &$1.90218 \pm5 \times 10^{-6}$ &Shanks, Wrench\\
1974 &$8\times10^{10}$ &$1.90216 04 \pm5 \times 10^{-7}$& Brent\\
1996 &$10^{14}$     &$1.90216 05778 \pm2.1 \times 10^{-9}$ &Nicely\\
2001 &$3\times10^{15}$ &$1.90216 05823 \pm 8 \times10^{-10}$ &Nicely\\
2007 &$7.7\times10^{15}$&$1.90216 05829 18\pm1.291 \times 10^{-9}$ &Nicely\\
2002 &$10^{16}$     &$1.90216 05831 04...$ &Sebah\\
\end{tabular}
\end{table}
ллллллллллллллл

Leaving only the primes fails to force convergence. If we pursue this thread, the most natural next step is to leave only the twin primes, that is, consecutive pairs of primes; it is customary (but not universal) to ignore 2 for this purpose and to count 5 twice, so the pairs are (3,5), (5,7), (11,13),..., (1019,1021),..., and these are incredibly sparse. In fact, it is not even known whether there is an infinite number of them and therefore if our series is infinite (this is called the Twin Primes Conjecture). Using only twin primes, all that is left of harmonic series is
$$
(\frac13+\frac15)+(\frac15+\frac17)+\cdots
$$
Do we achieve convergence now? Finally, the answer is yes, but no one is sure
to exactly what number; it is about $1.902,160,582,4\ldots$ and is known as Brun's
constant, after the Norwegian mathematician, Viggo Brun (1885-1978), who,
in 1919, established the convergence. Not much is known about it, although its
size is a strong indicator of just how sparse twin primes are. Thomas Nicely
provided the above estimate in 1994 and in the process uncovered the infamous
and much-publicized Intel Pentium division bug ('for a mathematician to get this
much publicity, he would normally have to shoot someone'), which made itself
apparent with the pair of twin primes 824,633,702,441 and 824,633,702,443.
His announcement to the world was by a now famous email, which began:
\begin{quote}
    It appears that there is a bug in the floating point unit (numeric
coprocessor) of many, and perhaps all, Pentium processors. In
short, the Pentium FPU is returning erroneous values for certain
division operations. For example, 1/824,633,702,441.0 is calculated
incorrectly (all digits beyond the eighth significant digit are
in error)\ldots
\end{quote}
On 17 January 1995 Intel announced a pre-tax charge of \S475 million against
earnings, as the total cost associated with the replacement of the flawed chips.

Incidentally, it would have been very convenient had the series diverged, as
that would have meant that there is an infinite number of twin primes and so
resolved the Twin Primes Conjecture. ( 5 is the only candidate for repetition by reasoning that any prime greater
than 3 must be of the form $6n\pm1$, any pair of twin primes must be 6n - 1 and
6n + 1 and therefore that a consecutive sequence of three is impossible beyond
3,5,7.) \cite[pp. 30, 31]{Havil}
лллллллллллллллллллллл

Brun's constant is defined to be the sum of the reciprocals of all twin primes
$$
B=\left(\frac{1}{3}+\frac{1}{5}\right)+\left(\frac{1}{5}+\frac{1}{7}\right)+\left(\frac{1}{11}+\frac{1}{13}\right)+\cdots
$$
If this series were divergent, then a proof of the twin prime conjecture would follow immediately. Brun proved, however, that the series is convergent and thus B is finite. His result demonstrates the scarcity of twin primes relative to all primes (whose reciprocal sum is divergent), but it does not shed any light on whether the number of twin primes is finite or infinite.

Selmer (1942), Fr\"{o}berg (1961), Bohman(1973) , Shanks $\&$ Wrench(1974), Brent(1975, 1976), Nicely(1995, 2001), Sebah, and others successively improved numerical estimates of B. The most recent calculations give
$$
B = 1.9021605831\ldots
$$
using large datasets of twin primes and assuming the truth of the extended twin prime conjecture. Let us elaborate on the latter issue. Under Hardy $\&$ Littlewood's hypothesis, the raw summation of twin prime reciprocals converges very slowly:
$$
\sum_{\substack{\text{twin} \\
p\leq n}}\frac1p-B=O\left(\frac{1}{\log n}\right),
$$
but the following extrapolation helps to accelerate the process
$$
\left(\sum_{\substack{\text{twin} \\
p\leq n}}\frac1p+\frac{4C_{twin}}{\log n}\right)-B=O\left(\frac{1}{\sqrt{n}\log n}\right)
$$
where $C_{twin} = 0.6601618158\ldots$ is the twin prime constant. Higher order extrapolations exist but do not present practical advantages as yet. In the midst of his computations, Nicely uncovered the infamous Intel Pentium error.
\cite[p. 133]{Finch}
ллллллллллллллл

%

ллллллллллллллл

%
лллллллллл

For every $x>1$, let
$$
\pi_2(x)= \#\{p\leq x:\ \ p+2\ \text{is also a prime}\}.
$$

Brun announced in 1919 that there exists an effectively computable integer $x_0$ such that, if $x \geq x_0$, then
$$
\pi_2(x)<100\frac{x}{\log^2x}
$$
The proof appeared in 1920.\\
In another paper of 1919, Brun had a weaker estimate for $\pi_2(x)$.

Based on heuristic considerations about the distribution of twin primes, B has been calculated by Kutrib and Richstein:
$$
B = 1.902160577783278....
$$

Brun also proved that for every $m \geq 1$ there exist $m$ successive primes which are not twin primes.\cite[p. 261]{Ribenboim}

лллллллллллл

To begin, in his famous paper of 1920, Brun showed that 2 may be written, in infinitely many ways, in the form $2=m-n$, where $m$, $n$ are 9-almost-primes\footnote{Let $k\geq1$. An integer $n =\prod_{i=1}^r p_i^{e_i}$ is called a k-almost-prime when $\sum_{i=1}^r e_i\leq k$.}. This was soon improved by Rademacher (1924), with 9 replaced by 7.

Later, R\'{e}nyi showed in 1947 that there exists $k\geq 1$ such that 2 may be written in infinitely many ways, in the form $2 = m - p$, where $p$ is a prime and $m$ is a k-almost-prime.

The best result to date, with sieve methods, is due to Chen (announced in 1966, published in 1973,1978); he proved that in R\'{e}nyi's result k may be taken equal to 2; so 2 = m - p, with m be 2-almost-prime and p prime, in infinitely many ways; this is very close to showing that there are infinitely many twin primes. A proof of Chen's theorem is given in the book of Halberstam and Richert. See also the simpler proof given by Ross (1975).

The sieve methods used for the study of twin primes are also appropriate for the investigation of Goldbach's conjecture.

Addendum on Polignac's Conjecture

The general Polignac conjecture can be, in part, treated like the twin-primes conjecture.

For every $k\geq1$ and $x > 1$, let $\pi_{2k}(x)$ denote the number of integers $n > 1$ such that $p_n\leq x$ and $p_{n+1}-p_n = 2k$.
With Brun's method, it may be shown that there exists a constant $C'_k > 0$ such that
$$
\pi_{2k}(x)<C'_k\frac{x}{\log^2x}
$$
\cite[pp. 264, 265]{Ribenboim}
ллллллллллллллллллллллллллллл


лллллллллллллл
\subsection{Gaps Between Primes}
Another result approaching the Twin Prime Conjecture is the recent work of Goldston, Pintz, and Yildirim, showing that
$$
\liminf_{n\rightarrow\infty}\frac{p_{n+1}-p_n}{\log p_n}=0
$$
where here, as through this paper, $p_n$ represents the $n$-th prime.\cite[p. 16]{Klyve}
лллллллллллллллллл

A subject that has attracted attention, but concerning which the known results leave much to be desired, is that of the behavior of $p_{n+1}-p_n$. As regards a universal upper bound for this difference, the first result was found by Hoheisel, who proved that there exists a constant $\alpha<1$, such that
$p_{n+1}-p_n = O(p_n^\alpha)$. The best result so far known is due to Ingham, who showed that this estimate holds for any $\alpha$ greater than 38/61. In both cases, what is actually proved is that
$$
\pi(x+x^\alpha)-\pi(x)\sim\frac{x^\alpha}{\log x},\quad\text{as}\ x\rightarrow\infty.
$$

In a crude sense one can say, in view of the prime number theorem, that the average of $p_{n+1}- p_n$ is $\log p_n$. Erd\H{o}s was the first to prove that there are infinitely many $n$ for which $p_{n+1} - p_n$ is appreciably greater than $\log p_n$, and Rankin proved that there are infinitely many $n$ for which
$$
p_{n+1}-p_n>c(\log p_n)\frac{(\log_2p_n)(\log_4p_n)}{(\log_3p_n)^2}
$$
where $\log_2x=\log\log x$ and so on, and $c$ is a positive constant. In the opposite direction, Bombieri and Davenport(1966) proved that there are infinitely many $n$ for which
$$
p_{n+1}-p_n<(0.46\cdots)\log p_n
$$
Of course, if the "prime twins" conjecture is true, there are infinitely many $n$ for which $p_{n+1}- p_n = 2$.

There is a somewhat paradoxical situation in connection with the limit points of the sequence
$$
\frac{p_{n+1}- p_n}{\log p_n}
$$
Erd\H{o}s,and Ricci (independently) have shown that the set of limit points has positive Lebesgue measure, and yet no number is known for which it can be asserted that it belongs to the set.\cite[pp. 173, 174.]{Davenport}
ллллллллллллллл

ллллллллллллллл
лллллллллллллллллл

It has been known for a long time that
$$
\limsup_{n\rightarrow\infty}\frac{p_{n+1}-p_n}{\log p_n}=\infty
$$
However, we have
$$
\liminf_{p\rightarrow\infty}\frac{p_{n+1}-p_n}{\log p_n}\leq0.248
$$
Is further improvement possible? If the twin prime conjecture is true, then the limit infimum is clearly 0.\cite[p. 131]{Finch}
лллллллллллллллл

\subsection{Hardy and Littlewood}
A statement equivalent to the Prime Number Theorem is that an integer near $x$ has a $1/ \log x$ ``probability'' of being prime. If we wish to know the probability that $p$ and $p+2$ are both prime, where $p$ is near $x$, we might simply multiply the probability that each is prime individually to get a probability of $1/ \log^2x$. We could then ``sum'' over all primes up to x to conjecture that
$$
\pi_2(x)\sim li_2(x):=\int_2^x\frac{1}{\log^2t}dt
$$
There is one crucial flaw with this argument, however, namely that $p$ and $p + 2$ being prime are not independent events. Consider, for example, that if $p > 2$ is prime, then $p$ is odd, and therefore $p + 2$ is also odd, giving it an immediate leg up on being prime. Further, in order to be prime, $p + 2$ must not be divisible by any odd prime $q$. A random integer has a $\ (1-1/q)$ chance of being Уnot divisibleФ by $q$, but if $q\nmid p$, then $p + 2$ must fall into one of $(q - 2)$ out of $(q - 1)$ remaining residue classes in order to be Уnot divisibleФ by $q$, and therefore we would expect that $p+2$ has a $(q-2)/(q-1)$ chance of being not divisible by $q$. If we assume that the primes are randomly distributed, we should then add a correction factor of
$$
\frac{p/(p-1)}{(p-1)/(p-2)}
$$
for each odd prime to the above approximation. Let us then amend our guess to
$$
\pi_2(x)\sim2\prod_{p>2}\frac{p(p-2)}{(p-1)^2}.li_2(x)=2\alpha.li_2(x)
$$
where $\alpha$ is called the Twin Prime Constant, and
$$
\alpha=\prod_{p>2}\left(1-\frac{1}{(p-1)^2}\right)\approx0.6601618158\ldots
$$

\begin{xconj}[Hardy and Littlewood]
For every integer $k > 0$ there are infinitely many prime pairs $p, p+2k$, and the number $\pi_{2k}(x)$ of such pairs less than $x$ is
$$
\pi_{2k}(x)\sim2\alpha\prod_{\substack{p>2 \\p|k}}\frac{p-1}{p-2}.li_2(x)
$$
\end{xconj}

de PolignacТs conjecture, of course, concerned \emph{consecutive} primes which differ by an even number, while Hardy and LittlewoodТs conjecture concerns any pair of primes which differ by an even number. However, it is not difficult to use the full version of the Hardy-Littlewood conjecture to show that the same asymptotic estimate given in Conjecture above holds for consecutive primes as well.

It is of some interest to note that the work of Hardy and Littlewood implies that $\pi_2(x)\sim\pi_4(x)$, and indeed there has been some attention given to the functions $\pi_{2k}(x)$ for small values of $k$. Those counted by $\pi_4(x)$ are often referred to as cousin primes.\cite[pp. 12, 13]{Klyve}\\
ллллллллллллллл

Striking theoretical progress has been achieved toward proving these conjectures(Goldbach's Conjecture and Twin Prime Conjecture), but insurmountable gaps remain. These formulas attempt to answer the following question:

Putting aside the existence issue, what is the distribution of primes satisfying various additional constraints? In essence, one desires asymptotic distributional formulas analogous to that in the Prime Number Theorem.

\textbf{Extended Twin Prime Conjecture.}
$$
\pi_2(n)\sim2C_{twin}\frac{n}{\log^2n}
$$
where
$$
C_{twin}=\prod_{p>2}\frac{p(p-2)}{(p-1)^2}=0.6601618158... =\frac12(1.3203236316...)
$$

Riesel(1985) discussed prime constellations, which generalize prime triplets and quadruplets, and demonstrated how one computes the corresponding Hardy-Littlewood constants. He emphasized the remarkable fact that, although we do not know the sequence of primes in its entirety, we can compute Hardy-Littlewood constants to any decimal accuracy due to a certain transformation in terms of Riemann's zeta function $\zeta(x)$.\cite[pp. 84-86]{Finch}\\
лллллллллллллллл


There are many who believe that Hardy's view is as valid now as it was then. Nevertheless, many dramatic advances have taken place in the intervening years.
\cite[p. 1]{Halberstam}
лллллллллл

A famous conjecture (1923) of Hardy and Littlewood is that $\pi_2(x)$ increases much like the function
$$
L_2(x)=2C\int_2^x\frac{du}{\log^2u}
$$
where $C = 0.661618158...$ is known as the twin-prime constant. The next table gives some idea how closely $\pi_2$ is approximated by $L_2(x)$.\cite[pp. 375, 376]{Burton}

\begin{table}[h!]
  \centering
  \caption{$\pi_2(x)$ approximation by $L_2(x)$}\label{constant}
\begin{tabular}{l l l l}\\\\
$x$ &$\pi_2(x)$& $L_2(x)-\pi_2(x)$\\
\hline
$10^3$&35&11\\
$10^4$&205&9\\
$10^5$&1,224&25\\
$10^6$&8,169&79\\
$10^7$&58,980&-226\\
$10^8$&440,312&56\\
$10^9$&3,424,506&802\\
$10^{10}$&27,412,679&-1262\\
$10^{11}$&224,376,048&-7183
\end{tabular}
\end{table}

лллллллллллллл

The estimate of $\pi_2(x)$ has been refined by a determination of the constant and of the size of the error. This was done, among others, by Bombieri and Davenport, in 1966. It is an application of sieve methods and its proof may be found, for example, in \cite{Halberstam}.

Here is the result:
$$
\pi_2(x)\leq8\prod_{p>2}\left(1-\frac{1}{(p-1)^2}\right)\frac{x}{\log^2x}\left\{1+O\left(\frac{\log\log x}{\log x}\right)\right\}
$$
The constant 8 has been improved to $68/9 + \varepsilon$ by Fouvry and Iwaniec (1983) and further by Bombieri, Friedlander, and Iwaniec (1986) to $7 + \varepsilon$. S. Lou has told to Ribenboim that he succeeded in obtaining a constant $6.26 + \varepsilon$, but he has not seen a paper with this result. Actually, it was conjectured by Hardy
and Littlewood that the factor should be 2 instead of 8.
\cite[pp. 261, 262]{Ribenboim}
ллллллллллллллл

Properties of $\pi(x)$ may be proved using the well-known relationship between the distribution of primes and the location of the zeros of the Riemann zeta function. Unfortunately, no similar relationship is known for twin primes, so very little is known about $\pi_2(x)$. It is not known whether there are infinitely many twin primes, and much less whether
$$
\pi_2(x)\sim L_2(x)
$$
as $x\rightarrow\infty$. However, empirical evidence suggests that this formula is true.\cite{Brent}
лллллллллллллллл
\subsection{Computations and Tables}

\begin{xconj}[de Polignac]
Every even number is the difference of two consecutive primes in infinitely many ways.
\end{xconj}

Taking the even number in the conjecture to be 2 immediately gives what
we now call the Twin Prime Conjecture.

de PolignacТs conjecture seems to have had very little impact when it
was published, and it was 30 years before the subject was revisited in the
literature. Glaisher inaugurated a project which continues to the present day
by enumerating the twin primes up to $10^5$.

Glaisher(1878) used published tables of primes to show that $\pi_2(10^5) = 1224$. In fact, he was interested in the density of twin primes, or as he called it, ``the rapidity of the decrease'' of the frequency of twin primes ``as we ascend higher in the series of numerals.'' Glaisher tabulated the number of twin primes in each of the intervals $[i\cdot10^6, i\cdot10^6 +100,000]$ for $i\in\{0, 1, 2, 6, 7, 8\}$, which was as far as his tables extended (why those tables didn't cover the range between 3,000,000 and 6,000,000 is not mentioned). Despite the large amount of data he compiled, Glaisher did not attempt to give an estimate for the density of the twin primes, claiming only that the ``number of prime-pairs\ldots [is] rather less than one-tenth of the number of primes in the same intervals.''He went on to claim that ``there can be little or no doubt that the number of prime-pairs is unlimited'' and that``it would be interesting, though probably not easy, to prove this.''(!)\cite[pp. 8, 9]{Klyve}\\\\
лллллллллллллл
\textbf{Lower bounds on $\pi_2(x)$}:

To this day we had no lower bound on $\pi_2(x)$ for large $x$, except for the finite number of twin primes that have been discovered. This is unsurprising in light of the fact that that we cannot even prove that there are infinitely many twin primes. Perhaps the closest we have come is in the work of Chen and others, who have studied the function
$$
\pi_{1,2}(x)=\#\{p\leq x:\ \Omega(p+2)\leq2\}
$$
where $\Omega(n)$ indicates the total number of prime factors of $n$, so that, for example, $24=2^3\cdot3^1$ and $\Omega(24) =3+1=4$. This function, $\pi_{1,2}(x)$, is said to count the number of twin almost primes. The number of these twin almost primes was proved by Chen to be infinite in 1973. Chen also gave the lower bound
$$
\pi_{1,2}(x)\geq0.335\cdot2\alpha li_2(x)
$$
for $x$ sufficiently large. This lower bound has been improved several times since, with the best result known given by Wu(2004) who showed
$$
\pi_{1,2}(x)\geq1.104\cdot2\alpha li_2(x)
$$
It is interesting to note that because of this result, it has now been demonstrated that the number of twin almost primes is greater than the conjectured number of twin primes. \\\\
лллллллллллллллл
\textbf{Computing $\pi_2(x)$}.

Following the work of Hardy and Littlewood, the story of the attempts to compute $\pi_2(x)$ to ever-higher values of $x$ is one of testing a conjecture that everyone believes. Notable among recent calculations are those of Nicely and Sebah. Each successive computation has served to give only greater confirmation for the strong form of the Twin Prime Conjecture. Indeed Sebah's computation up to $10^{16}$ (the current record holder) gives that
$$
\pi_2(x)=10,304,195,697,298
$$
while an estimate using Hardy-Littlewood would leave us to expect that
$$
\pi_2(x)\approx10,304,192,554,496
$$
The uncanny precision of the strong form of the Twin Prime Conjecture is easily seen, and it now seems all but impossible to believe that this conjecture could be false.\cite[p. 17]{Klyve}\\\\
лллллллллллллл
\textbf{Conjectures involving two different kinds of prime triplets:}
$$
P_n(p,p+2,p+6)\sim P_n(p,p+4,p+6)\sim D\frac{n}{\log^3n}
$$
where
$$
D=\frac92\prod_{p>3}\frac{p^2(p-3)}{(p-1)^3}=2.8582485957\ldots
$$
\textbf{Conjectures involving two different kinds of prime quadruplets:}
$$
P_n(p,p+2,p+6,p+8)\sim \frac12P_n(p,p+4,p+6,p+10)\sim E\frac{n}{\log^4n}
$$
where
$$
E=\frac{27}{2}\prod_{p>3}\frac{p^3(p-4)}{(p-1)^4}=4.1511808632\ldots
$$
\textbf{Conjecture involving primes of the form $m^2+1$:}

If $Q_n$ is defined to be the number of primes $p\leq n$ satisfying $p = m^2 + 1$ for some integer $m$, then
$$
Q_n\sim2C_{quad}\frac{\sqrt{n}}{\log n}
$$
where
$$
C_{quad}=\frac12\prod_{p>2}\left(1-\frac{(-1)^{\frac{p-1}{2}}}{p-1}\right)=0.6864067314\ldots = \frac12(1.3728134628\ldots).
$$
\textbf{Extended Goldbach Conjecture:}

If $R_n$ is defined to be the number of representations of an even integer $n$ as a sum of two primes (order counts), then
$$
R_n\sim2C_{twin}\cdot\frac{n}{\log^2n}\prod_{\substack{p>2\\p|n}}\frac{p-1}{p-2}
$$
It is intriguing that both the Extended Twin Prime Conjecture and the Extended Goldbach Conjecture involve the same constant $C_{twin}$. It is often said that the Goldbach conjecture is ``conjugate'' to the Twin prime conjecture.\cite[p. 85]{Finch}
ллллллллллллллл

\begin{xthm}[Ramar\'{e}. 1995,1996]
Every even integer can be expressed as a sum of six or fewer primes (in other words, Schnirelmann's number is $\leq6$).
\end{xthm}
\begin{xthm}[Chen. 1973]
Every sufficiently large even integer can be expressed as a sum of a prime and a positive integer having $\leq2$ prime factors.
\end{xthm}
In fact, Chen proved the asymptotic inequality
$$
\widetilde{R}_n>0.67\cdot\frac{n}{\log^2n}\cdot\prod_{p>2}\left(1-\frac{1}{(p-1)^2}\right)\cdot\prod_{\substack{p>2\\p|n}}\frac{p-1}{p-2}
$$
where $\widetilde{R}_n$ is the number of corresponding representations. Chen also proved that there are infinitely many primes p such that p+2 is an almost prime, a weakening of the twin prime conjecture, and the same coefficient 0.67 appears.\cite[p. 88]{Finch}
лллллллллллллллл

Many examples of immense twins are known. The largest twins to date, each 51090 digits long,
$$
33218925\cdot 2^{169690} \pm 1
$$
were discovered in 2002.\cite[p. 50]{Burton}
ллллллллллллллл

Table \ref{growth} gives a feeling for the growth of $\pi_2(x)$. I reproduce part of the calculations by Brent (1975,1976).
\begin{table}[h!]
  \centering
  \caption{growth of $\pi_2(x)$}\label{growth}
\begin{tabular}{l  l}\\\\
$x$ & $\pi_2(x)$\\
\hline
$10^3$&35 \\
$10^4$&205\\
$10^5$&1224\\
$10^6$&8169\\
$10^7$&58980\\
$10^8$&440312 \\
$10^9$&3424506 \\
$10^{10}$&27412679 \\
$10^{11}$&224376048 \\
$10^{12}$&1870585220 \\
$10^{13}$&15834664872 \\
$10^{14}$&135780321665 \\
$1.37\cdot10^{14}$&182312485795 \\
\end{tabular}
\end{table}

The largest exact value for the number of twin primes below a given ("round number") limit is $\pi_2(1.37\cdot10^{14})$. The calculation to this bound is due to T.R. Nicely (communicated by letter of October 10, 1995). On the other hand M. Kutrib and D. Richstein reached independently the bound $10^{14}$ in September, 1995. The results of these calculations agreed completely.

These extensive calculations were preceded by - for their time - just as extensive calculations(Table \ref{p2 person}).
\begin{table}[h!]
  \centering
  \caption{constant}\label{p2 person}
\begin{tabular}{l  l}\\\\
\hline
$\pi_2(10^5)$& by J.W.L. Glaisher (1878), \\
$\pi_2(10^{6})$& by G.A. Streatfeild (1923), \\
$\pi_2(37\cdot10^{6})$& by D.H. Lehmer (1957), \\
$\pi_2(104.3\cdot10^{6})$& by G. Armendiny and F. Gruenberger (1961), \\
$\pi_2(2\cdot10^{8})$& by S. Weintraub (1973), \\
$\pi_2(2\cdot10^{9})$& by J. Bohman (1973), \\
$\pi_2(10^{11})$& by R.P. Brent (1995).
\end{tabular}
\end{table}

There is keen competition to produce the largest pair of twin primes.

On October 9, 1995, Dubner discovered the largest known pair of twin primes (with 5129 digits): (p, p + 2) where
$$
p = 570918348 \cdot 10^{5120} - 1.
$$
It took only one day (with 2 crunchers) (the expected time would be 150 times longer! What luck!).

Earlier, in 1994, K.-H. Indlekofer and A. J\'{a}nai, found the pair $697,053,813\cdot 2^{16352}\pm1$ (with 4932 digits). In 1995 the amateur mathematician T. Forbes found the third largest known pair $6,797,727\cdot2^{15328}\pm 1$ (with 4622 digits). These surpass the pair found in 1993 by Dubner to be $1,692,923,232 \cdot 10^{4020}\pm1$. These numbers have 4030 digits.

A list of twin primes with at least 1000 digits (Titanic primes) was originally compiled by S. Yates. After he passed away, this work has been continued by C. Caldwell. This list contained 28 pairs of twin primes; it gives the discoverer, year of discovery, and number of digits of each prime.

The eight largest entries in the list are shown in Table \ref{largest}.
\begin{savenotes}
\begin{table}[h!]
  \centering
  \caption{constant}\label{largest}
\begin{tabular}{l l l l}\\\\
Twin Primes&Discoverer&of digits&Year\\
\hline
$570918348 \cdot 10^{5120}\pm1$&Dubner &5129 &1995 \\
$697053813\cdot 2^{16352}\pm1$&Indlekofer and Jсrai&4932 &1994 \\
$6797727\cdot 2^{15328}\pm1 $ &Forbes&4622 &1995 \\
$1692923232\cdot10^{4020}\pm1$&Dubner &4030 &1993 \\
$4655478828\cdot 10^{3429}\pm1$&Dubner &3439 &1993 \\
$1706595\cdot2^{11235}$&Just three from Amdahl\footnote{ B. Parady, J. Smith, S. Zarantonello. }&3389 &1989 \\
$459\cdot2^{8529}\pm1$&Dubner &2571 &1993\\
$1171452282\cdot 10^{2490}\pm1$&Dubner &2500 &1991
\end{tabular}
\end{table}
\end{savenotes}
Besides the names already quoted, the list also includes primes discovered by Keller, and Atkin and Rickert.\cite[pp. 262, 263]{Ribenboim}
лллллллллллллллллллллл
%
%
%

\bibliographystyle{plain}

\begin{thebibliography}{99}
\bibitem{Andrews} George E. Andrews, Sieves for theorems of Euler, Rogers, and Ramanujan. The theory of arithmetic functions (Proc. Conf., Western Michigan Univ., Kalamazoo, Mich., 1971), pp. 1Ц20. Lecture Notes in Math., Vol. 251, Springer, Berlin, 1972.
\bibitem{Bateman} Paul T. Bateman, Harold G. Diamond, Analytic number theory.
An introductory course. World Scientific Publishing Co. Pte. Ltd., Hackensack, NJ, 2004.
\bibitem{Brent} Richard P. Brent, Irregularities in the Distribution of Primes and Twin Primes, MATHEMATICS OF COMPUTATION, VOLUME 29, NUMBER 129, JANUARY 1975, PAGES 43-56
\bibitem{Burton} Burton, David M. Elementary number theory. McGraw-Hill, 2007.
\bibitem{Chandrasekharan}  K. Chandrasekharan, Arithmetical functions. Die Grundlehren der mathematischen Wissenschaften, Band 167 Springer-Verlag, New York-Berlin 1970
\bibitem{Davenport} Harold Davenport, Multiplicative number theory. Second edition. Revised by Hugh L. Montgomery. Springer-Verlag, New York-Berlin, 1980.
\bibitem{Erdos1}Paul Erdos, Janos Suranyi, Topics in the Theory of Numbers,2003 Springer
\bibitem{Finch} Steven R. Finch, Mathematical constants. Encyclopedia of Mathematics and its Applications, 94. Cambridge University Press, Cambridge, 2003.
\bibitem{Halberstam} Halberstam, H.; Richert, H.-E. Sieve methods. London Mathematical Society Monographs, No. 4. Academic Press [A subsidiary of Harcourt Brace Jovanovich, Publishers], London-New York, 1974.
\bibitem{Havil}Havil, Julian, Gamma, Exploring Euler's constant. With a foreword by Freeman Dyson. Reprint of the 2003 edition [MR1968276]. Princeton Science Library. Princeton University Press, Princeton, NJ, 2009
\bibitem{Iwaniec} H. Iwaniec, Almost-primes represented by quadratic polynomials. Invent. Math. 47 (1978), no. 2, 171Ц188.
\bibitem{Klyve} Dominic Klyve, Explicit Bounds on Twin Primes and BrunТs Constant, A Thesis
Submitted to the Faculty in partial fulfillment of the requirements for the degree of Doctor of Philosophy in Mathematics
\bibitem{Landau} Edmund Landau, Elementary number theory. Translated by J. E. Goodman. Chelsea Publishing Co., New York, N.Y., 1958.
\bibitem{Landau c} E. Landau, Collected works, Thales Verlag, Essen, 1987.
\bibitem{Leveque} Leveque, Fundamentals of Number Theory
\bibitem{Nicely} Thomas R. Nicely, Enumeration to $10^{14}$ of the Twin Primes and Brun's
Constant, Virginia Journal of Science, Volume 46, Number 3, Fall 1996
\bibitem{Niven} Niven, Ivan; Zuckerman, Herbert S.; Montgomery, Hugh L. An introduction to the theory of numbers. Fifth edition. John Wiley $\&$ Sons, Inc., New York, 1991.
\bibitem{Ribenboim} Ribenboim, Paulo, The new book of prime number records. Springer-Verlag, New York, 1996.
\end{thebibliography}

\end{document}